\documentclass[10pt]{amsart}
\usepackage{amssymb}

\newtheorem{theorem}{Theorem}[section]
\newtheorem{proposition}[theorem]{Proposition}
\newtheorem{lemma}[theorem]{Lemma}
\newtheorem{corollary}[theorem]{Corollary}
\theoremstyle{definition}
\newtheorem{definition}[theorem]{Definition}

\newtheorem{remark}[theorem]{Remark}

\newtheorem*{thmA}{Theorem A}
\newtheorem*{thmB}{Theorem B}

\newcommand{\Z}{\mathbb Z}
\newcommand{\C}{\mathbb C}
\newcommand{\R}{\mathbb R}
\newcommand{\PP}{\mathbb P}
\newcommand{\cL}{\mathcal L}
\newcommand{\cI}{\mathcal I}
\newcommand{\cB}{\mathcal B}
\newcommand{\Mod}{\operatorname{Mod}}
\newcommand{\hX}{\widehat X}
\newcommand{\tL}{\widetilde L}
\newcommand{\Dred}{D_{\mathrm{red}}}

\title[Liftable braid factorizations]{Liftable braid factorizations from line arrangements:\\ Hironaka's relations and a Zariski pair}
\author{Atsuhide Mori}
\address{Department of Mathematics, Osaka Dental University}
\email{mori-a@cc.osaka-dent.ac.jp}

\begin{document}
\raggedbottom

\begin{abstract}
We realize Hironaka's generalized lantern relations by lifting positive
braid factorizations. For a nonpencil real line arrangement with no parallel
lines, we construct a family of branched coverings on the regular fibers of a
blown-up pencil of spheres. The motion of the branch points gives a
factorization of a power of the full twist into least
liftable powers of half-twists. The
Artal--Carmona--Cogolludo--Marco Zariski pair
gives two such factorizations of the same braid for a single simple
branched covering of the disk. Their lifts are Hurwitz inequivalent even
up to simultaneous conjugation by page homeomorphisms. The corresponding
Milnor fibers fill the same planar contact three-manifold and have the
same incidence matrix. They are homeomorphic and lie in distinct strong
symplectic deformation classes.
\end{abstract}

\maketitle

\section{Introduction}\label{sec:intro}
The braid monodromy of a complex line arrangement determines its embedded
topology \cite[Theorem~3.8]{Cogolludo}. Artal, Carmona, Cogolludo, and Marco
used braid monodromy to distinguish two complexified real arrangements
of eleven lines with the same combinatorics \cite{ACCM}.

Hironaka used the braid monodromy of a real line arrangement to express a
boundary multitwist on a planar surface as a product of positive Dehn
twists, one for each intersection point \cite{Hironaka}. We recover this
relation by lifting powers of half-twists through a simple branched
covering of the disk (Section~\ref{sec:layers}).

\begin{thmA}[Theorem~\ref{thm:arrangement}]
Let $\cL$ be a nonpencil real line arrangement with no parallel lines, in
the coordinates of Section~\ref{sec:arrangement}. Let $k$ be a compatible
positive half-integer (Definition~\ref{def:admissiblek}).
Then there is a simple branched covering of the disk with planar page for which 
the $2k$-th power of the full twist admits a positive factorization into
least liftable powers of half-twists.
The Dehn twist lifts recover Hironaka's factorization for $\cL$.
The factors with identity lift correspond to the nodes and cusps
of the branch curve of a fiberwise covering of the Hirzebruch surface $F_{2k}$.
\end{thmA}

Blowing up the intersections separates the lines into disjoint sections
of a pencil of spheres. We use these sections as poles of a map on each
fiber. The pole orders determine the boundary covering, and compatibility
allows the numerator to stay nonzero at the poles
(Section~\ref{sec:construction}). Each blow-up on a section changes its
boundary framing by one twist (Section~\ref{sec:framings}).

\begin{thmB}[Theorem~\ref{thm:zariskiplays}]
For the Zariski pair of Artal--Carmona--Cogolludo--Marco \cite{ACCM},
apply Theorem~A with $k=30$, the least value compatible with both
arrangements. This gives two positive factorizations of the $2k$-th power
of the full twist, with $2k=60$, for a single simple branched covering of
the disk. Their lifts are not Hurwitz equivalent, even up to simultaneous
conjugation by homeomorphisms of the page.
\end{thmB}

Collapsing the inner boundary circles to marked points recovers the braid
monodromy of each arrangement (Lemma~\ref{lem:collapsed-monodromy}).
The two braid monodromies distinguish the embeddings of the arrangements
and therefore distinguish the lifted factorizations.
Baykur and Van Horn-Morris give $4$-braids with infinitely many quasipositive
factorizations for a fixed double cover by a torus with two boundary
components. Homology distinguishes the total spaces of their Lefschetz
fibrations \cite[Theorem~10]{BVHM}.

The lifted factorizations in Theorem~B give Milnor fibers of rational surface
singularities with the same contact link and incidence matrix.
The boundary identification extends to a homeomorphism of the fillings
(Proposition~\ref{prop:milnor} and Corollary~\ref{cor:fillings}).
The inequivalence of the planar Lefschetz fibrations puts the fillings in
distinct strong symplectic deformation classes. Baykur, N\'emethi, and
Plamenevskaya observe that
different decorated disk arrangements can have the same incidence matrix
\cite[\S5.2]{BNP}. In our example, the contact link and incidence matrix
together fail to determine the strong symplectic deformation class.
Whether the underlying smooth four-manifolds are diffeomorphic remains
open.

\section{Liftable factorizations}\label{sec:layers}
Let $p\colon P\to D$ be a degree-$d$ simple branched covering of the disk
with branch set $\Delta$ of size $n$. The page $P\cong\Sigma_{0,r}$ is a
sphere with $r$ boundary components. Let $\sigma$ be the monodromy
permutation along $\partial D$. Braids on $\Delta$ that lift to
homeomorphisms fixing $\partial P$ pointwise form a subgroup $L_p$.
Taking the mapping class of the lift defines
$\ell_p\colon L_p\to\Mod(P\ \mathrm{rel}\ \partial P)$.
A \emph{chord} is an embedded arc $a\subset\operatorname{int}D$ joining two
distinct points of $\Delta$ and otherwise disjoint from $\Delta$.
Choose a regular base point in a small disk around $a$ and label the
sheets over it. Take arcs from that point to the endpoints, disjoint
except at the base point. Read the local monodromies as transpositions
$t,t'$ with respect to this common labeling.
Let $\sigma_a$ be the positive half-twist along $a$.
Its least positive power that lifts to $P$ fixing $\partial P$ pointwise is $h_a=\sigma_a^{e(a)}$,
where the \emph{exponent} $e(a)$ is $1$, $2$ or $3$ according as $t=t'$, $t$ and $t'$ are disjoint, or
$t$ and $t'$ share exactly one letter \cite{MP}. 
The preimage of a small disk around $a$ is a union of disks, 
except that for $e(a)=1$ the component over the two sheets exchanged by $t$ is an annulus. 
Hence $h_a$ lifts to the Dehn twist about the core of that annulus if $e(a)=1$, 
and to the identity if $e(a)\in\{2,3\}$. 
We write $\tau_c$ for a positive Dehn twist about a curve $c$.
A \emph{liftable factorization} of a braid $\beta\in L_p$ is a sequence
$w=(h_{a_1},\ldots,h_{a_m})$ of least liftable powers of half-twists whose product is $\beta$.
Its \emph{lifted factorization} is the ordered list of the Dehn twists
$\ell_p(h_{a_j})$ with $e(a_j)=1$, and its product is $\ell_p(\beta)$.
Two liftable factorizations are \emph{equivalent} if they are related by Hurwitz moves and 
simultaneous conjugation by elements of $L_p$.
For factorizations of a fixed braid $\beta$, 
simultaneous conjugation is restricted to elements of $C_{L_p}(\beta)$, 
which is all of $L_p$ when $\beta$ is central. 
Lifting respects Hurwitz moves and simultaneous conjugation, 
so equivalent liftable factorizations have equivalent lifted factorizations.

Write $\tau_{\partial D}$ for the full twist, represented by a Dehn twist
in a collar of $\partial D$.
\begin{lemma}\label{lem:fulltwist}
Let $m\ge1$ and suppose that every cycle length of $\sigma$ divides $m$.
Then $\tau_{\partial D}^{\,m}\in L_p$ and
\[
 \ell_p\bigl(\tau_{\partial D}^{\,m}\bigr)=\prod_{c}\tau_{\partial_c}^{\,m/|c|},
\]
where the product is over the cycles $c$ of $\sigma$, $|c|$ is the length of
$c$, and $\partial_c$ is the corresponding boundary circle of $P$.
\end{lemma}

\begin{proof}
Represent $\tau_{\partial D}^m$ by a twist in a collar $K$ of $\partial D$.
Each annulus $K_c$ over $K$ covers it with degree $|c|$, so the twist lifts
on $K_c$ through angle $2\pi m/|c|$. Extend these twists by the identity
off $p^{-1}(K)$. Since each $|c|$ divides $m$, the lift fixes $\partial P$
pointwise and represents $\prod_c\tau_{\partial_c}^{m/|c|}$.
\end{proof}

\section{Hironaka's relations on branched covers}\label{sec:arrangement}

Let $\cL=\{L_1,\ldots,L_s\}$ be a real line arrangement in $\R^2\subset\C^2$
with $s\ge3$, distinct slopes, no vertical line, and intersection points
$\cI$ with distinct $x$-coordinates, as in \cite[Theorem~1.2]{Hironaka}.  We
assume $\cL$ is not a pencil, so each $L_i$ contains at least two points of
$\cI$. Write $\mu_i$ for this number.
For $p\in\cI$, put $S_p=\{i:p\in L_i\}$.

Take a vertical reference fiber to the right of all the intersections.
Its page is a disk with small disks removed around its intersections
with the lines. The outer boundary is $\partial_0$, and the hole for
$L_i$ is $\partial_i$. In a fiber just to the right of $x(p)$, draw a
circle around the holes for the lines through $p$. Transport it back
to the reference page along Hironaka's path through the upper
half-plane to obtain $a_p$ \cite[Lemma~2.1]{Hironaka}.
Hironaka's relation is
\begin{equation}\label{eq:hironaka}
 \tau_{\partial_0}\prod_{i=1}^s\tau_{\partial_i}^{\,\mu_i-1}
 =\tau_{a_{p_{s'}}}\cdots\tau_{a_{p_1}},
\end{equation}
where $s'=\#\cI$ and the points $p_1,\ldots,p_{s'}$ are numbered in decreasing
order of $x$-coordinate \cite[Theorem~1.2]{Hironaka}.

\begin{definition}\label{def:admissiblek}
A number $k\in\frac12\Z_{>0}$ is \emph{compatible} with $\cL$ if
$(\mu_i-1)\mid 2k$ for all $i$.  For such $k$ put
\[
 \lambda_0=2k,\qquad \lambda_i=\frac{2k}{\mu_i-1}\ (1\le i\le s),\qquad
 d=\sum_{c=0}^s\lambda_c,\qquad n=d+s-1 .
\]
\end{definition}

By Riemann--Hurwitz, a connected simple branched covering of degree $d$ whose
boundary monodromy has cycle type $(\lambda_0,\ldots,\lambda_s)$ has a planar
page if and only if it has $n=d+s-1$ branch points, and such coverings exist.
We denote this covering by $p_{\cL,k}$. It is unique up to Hurwitz moves
and conjugation \cite[Theorem~1]{Kluitmann}.

\begin{theorem}\label{thm:arrangement}
Let $k$ be compatible with $\cL$.  Then $\tau_{\partial D}^{\,2k}$ lies in
$L_{p_{\cL,k}}$ and has a liftable factorization $w$ such that
\begin{enumerate}
\item the factors of exponent one correspond to the points of $\cI$, and, for
 a suitable identification of the page with Hironaka's surface under which
 the boundary circle over a cycle of length $\lambda_c$ corresponds to
 $\partial_c$ $(0\le c\le s)$, the lifted factorization of $w$ is
 Hironaka's factorization \eqref{eq:hironaka};
\item the other factors have exponent two or three.
\end{enumerate}
\end{theorem}

\begin{remark}\label{rem:lantern}
The boundary monodromy of $p_{\cL,k}$ is trivial exactly when all
$\lambda_c=1$, that is, when $2k=1$ and $\mu_i=2$ for all $i$.  This happens
only for three lines in general position.  In that case $d=4$ and $n=6$,
and the liftable factorization $w$ of the full twist lifts to the lantern relation.
\end{remark}

\section{The construction}\label{sec:construction}

At an intersection of the arrangement, the source fiber splits into two
spheres, and the branch points near their node make the half-twist that
lifts to a vanishing cycle twist.

\subsection{The blown-up pencil}
The projection $\rho(x,y)=x$ is the pencil of lines through
$O=[0:1:0]$. Blowing up $O$ gives a sphere bundle over $\PP^1$. Blowing up
the intersection points of $\cL$ gives $\pi\colon\hX\to\PP^1$, with a
reducible fiber $E_p\cup\widetilde F_p$ over $x(p)$ for each $p\in\cI$.
Here $E_p$ is the exceptional sphere and $\widetilde F_p$ is the proper
transform of the vertical line through $p$. All other fibers are spheres.
The blow-ups separate the lines into disjoint sections $S_i=\tL_i$,
and $S_0=E_O$ is a further section. Removing a small disk around each
section point on a smooth fiber gives Hironaka's planar page
\cite[\S2.2]{Hironaka}. Write $F$ for the fiber class and
$e_c=-S_c^2$. The self-intersections are $-e_0=-1$ and
$-e_i=1-\mu_i$, since a blow-up on a line lowers its self-intersection
by one.

\subsection{Prescribing the poles}
On each fiber we prescribe a pole of order $\lambda_c$ at the point on
$S_c$. In particular, a small circle about this pole covers a circle
near infinity with degree $\lambda_c$. The numerator must stay nonzero
along the pole sections, so its line bundle must have degree zero there.
Tensoring with a line bundle pulled back from the base adjusts these
degrees without changing the degrees on the fibers. We use
\[
 D_\infty=\sum_{c=0}^s\lambda_c S_c,\qquad D_0=D_\infty+2kF .
\]
The first divisor records the poles with their multiplicities. We
represent the denominator by a section $s_\infty$ with zero divisor
$D_\infty$ and choose the numerator $s_0$ as a section of $\mathcal O(D_0)$.
The ratio $s_0/s_\infty$ takes values in the pullback of
$\mathcal O_{\PP^1}(2k)$, so its target after adding infinity to each
fiber is the Hirzebruch surface $F_{2k}$.

Each section meets a fiber once, and the sections are pairwise disjoint.
Thus the numerator bundle has degree $D_0\cdot S_c=2k-\lambda_c e_c=0$
on $S_c$, by compatibility.

\begin{lemma}\label{lem:linsys}
A general numerator $s_0$ is nonzero on all the pole sections. The
restriction map for numerators is onto on each fiber, including the
reducible fibers.
\end{lemma}

\begin{proof}
Put $\Dred=\sum_c S_c$. The restriction sequences reduce both
assertions to the vanishing of $H^1$ for $D_0-\Dred$ and $D_0-F$.
For a divisor $G+jF$ with $G=\sum_c m_cS_c$ and $m_c\ge0$, remove one
copy of a section at a time. If $A$ is the divisor at the current
step and $S_c$ is the section to be removed, use
\[
0\longrightarrow\mathcal O(A-S_c)\longrightarrow\mathcal O(A)
\longrightarrow\mathcal O_{S_c}(A)\longrightarrow0 .
\]
The last term has degree at least $j-m_ce_c$, its degree before
any sections are removed. If $j-m_ce_c\ge-1$ for every $c$, the
restrictions have vanishing $H^1$. The remaining divisor $jF$ has
$H^1(\hX,\mathcal O(jF))=H^1(\PP^1,\mathcal O(j))=0$ for $j\ge-1$,
since the fibers are connected trees of spheres. The restriction
sequences for the removed sections give the desired vanishing.
For $D_0-\Dred$ the initial degrees are $e_c$, and for $D_0-F$ they
are $-1$. Finally, the numerator bundle has degree zero on each $S_c$
and is therefore trivial there, so a general numerator
restricts to a nonzero constant in a trivialization on each section.
\end{proof}

\subsection{Branch points and local braids}
The numerator and denominator have no common zeros. They define the
fiberwise map
\[
 \Phi=[s_\infty:s_0]\colon\hX\longrightarrow F_{2k}
\]
over $\PP^1$. The section at infinity is the negative section $E_{2k}$,
with $\Phi^*E_{2k}=D_\infty$. A smooth fiber maps with degree
$d=\sum_c\lambda_c$. In the fiber over $x(p)$, the poles indexed by
$i\in S_p$ lie on $E_p$ and the others lie on $\widetilde F_p$.
The degrees on the two components are therefore $\sum_{i\in S_p}\lambda_i$
and $d-\sum_{i\in S_p}\lambda_i$, both positive. Since every fiber
component maps with positive degree, $\Phi$ is finite. The finite critical values
over smooth source fibers sweep out a curve. Its closure in $F_{2k}$
is the \emph{residual branch curve} $\cB$.

\begin{lemma}\label{lem:generic}
For general $s_0$, the curve $\cB$ is disjoint from $E_{2k}$.
Away from the fibers over $x(\cI)$, its only singularities are ordinary
nodes and cusps.  The two local transpositions are disjoint at a node and
share exactly one letter at a cusp.  Over each $x(p)$, its only
nontransverse intersection with the fiber is a simple vertical tangency
at the image of the node of $E_p\cup\widetilde F_p$, with equal local
transpositions.
\end{lemma}

\begin{proof}
\emph{Poles.} In a target coordinate vanishing at infinity, the local
map at a pole of order $\lambda$ is $z^\lambda$ times a unit.
Consequently no other critical points accumulate on a pole section,
and $\cB$ misses $E_{2k}$.

\emph{Smooth fibers.} There are two relevant events. The values of two
distinct simple critical points can cross. If the crossing is transverse,
the branch curve has a node. The two ramification points lie in disjoint
degree-two local covers of the same disk about the common value. They
remain disjoint over nearby fibers, so the transpositions are disjoint.
Alternatively, two simple critical points can coalesce. The transverse
one-parameter model is $z^3+tz$, whose critical values trace a cusp.
The two transpositions act transitively on the three sheets of this
connected local cover, so they share exactly one letter.

To justify genericity, write $f=g/h$ on a smooth fiber, with $h$ fixed
and $g$ the restriction of the numerator. A simple critical point
satisfies $f'=0$ and $f''\ne0$. Lemma~\ref{lem:linsys} and polynomial
interpolation on $\PP^1$ let us prescribe jets of $f$ away from the
poles of total length at most $d+1$. For $j$ distinct points on a
varying fiber there are $1+j$ complex variables, including the base
coordinate. The following table gives the number $c$ of independent
conditions and the jet length needed for their independence.
\[
\begin{array}{l|ccc}
\text{Event} & j & c & \text{Jet length}\\ \hline
f'=f''=f'''=0 & 1 & 3 & 4\\
\text{A nonsimple critical point and another with the same value} & 2 & 4 & 5\\
\text{Three critical points with the same value} & 3 & 5 & 6\\
\text{Two simple critical points with the same value} & 2 & 3 & 4\\
f'=f''=0,\quad f'''\ne0 & 1 & 2 & 3
\end{array}
\]
For example, the fourth row imposes two derivative conditions and one
equality of values. Independence in the numerator parameters makes the
joint solution space smooth. A regular value of its projection onto
the numerator space gives a transverse solution set of dimension
$1+j-c$, and generic smoothness gives such regular values for general
numerators. Thus the first three rows, where $c>1+j$, are excluded when
the jets can be prescribed. Since $d\ge4$, the only exceptions to the
jet bound are the second and third rows for $d=4$. Degree excludes
these events, since the total multiplicities over the common value
would be at least $3+2$ and $2+2+2$, respectively.

The last two rows give finitely many events, since the solution sets are
zero-dimensional and algebraic. The Jacobian in the base and point
coordinates is invertible. For the fourth row this says that the
difference of the two critical values has a simple zero. For the fifth
row it gives $f_{xz}\ne0$ at the double critical point, so the family
is locally the unfolding $z^3+tz$. Away from these events,
simple critical values vary holomorphically with the base, so there are
no vertical tangencies over smooth fibers. Events on the same fiber
have distinct critical values and give commuting braids supported in
disjoint disks.

\emph{Reducible fibers.} The node represents the vertical direction
at $p$. The lines through $p$ meet the exceptional sphere in their
nonvertical directions, so the node lies off the pole sections.
Lemma~\ref{lem:linsys} lets us prescribe the
numerator on the two spheres subject to agreement at the node. On
each sphere, the same interpolation with the base fixed gives simple
critical points with distinct finite critical values. We can also make
the map unramified at the node and keep its value away from the finite
critical values. On one component, write the map as $g/h$ and let $c$
be its value at the node. The base twist $2kF$ restricts trivially to
this component, so $g$ and $h$ are sections of the same line bundle
after trivializing the base twist. Replacing $g$ by
$\lambda g+(1-\lambda)ch$ rescales the map about $c$ and separates its
critical values from those on the other component for general
$\lambda\in\C^\times$. On each reducible fiber, these conditions define
a nonempty Zariski-open set of restrictions. By Lemma~\ref{lem:linsys},
their inverse images are nonempty Zariski-open sets of numerators.
There are finitely many reducible fibers, so a general $s_0$ satisfies
all these conditions.

Near the node, write $t=uv$ and $w=c+au+bv+\cdots$, with $a,b\ne0$.
Use $(w,v)$ as coordinates in the source. As a function of $v$ at fixed
$w$, $t$ has a nondegenerate critical point near the node, since its
second derivative there is $-2b/a$. The resulting local degree-two map
has branch curve
\[
 t=\frac{(w-c)^2}{4ab}+O\bigl((w-c)^3\bigr).
\]
This gives a simple vertical tangency. The two branch points on a nearby
fiber belong to this local degree-two cover, so their transpositions agree.
\end{proof}

\subsection{Proof of Theorem~\ref{thm:arrangement}}
Fix a regular real base point $x_0$ to the right of $x(\cI)$, as in \cite{Hironaka}, and a
disk $D$ in the fiber of $F_{2k}$ over $x_0$ containing
$\cB\cap\Phi(\pi^{-1}(x_0))$ and missing $E_{2k}$.  Over $D$ the map $\Phi$
restricts to a simple branched covering of degree $d$ whose page is obtained
by removing neighborhoods of the poles from the fiber.  The pole orders
$\lambda_0,\ldots,\lambda_s$ give the cycle type of the boundary monodromy.
The page is $\Sigma_{0,s+1}$, so Riemann--Hurwitz gives $n=d+s-1$ branch
points.  This identifies the covering with $p_{\cL,k}$.

The covering extends over each loop based at $x_0$. Its lifted monodromy
fixes the boundary pointwise by the computation in the proof of
Lemma~\ref{lem:framings}, so the braid monodromy takes values in $L_p$.
Lemma~\ref{lem:generic} identifies the local factors as the least liftable
powers of exponents one at the tangencies over $x(p)$, two at nodes, and
three at cusps.
The factor over $x(p)$ lifts to the twist about the vanishing cycle $a_p$
of the corresponding singular fiber of $\pi$.  Hironaka's choice of loops
gives the prescribed order and hence the liftable factorization $w$.
For general $s_0$, the curve $\cB$ meets the fiber over $x=\infty$
transversely at smooth points, by the surjectivity in
Lemma~\ref{lem:linsys}. In the fiber coordinate $w'=w/x^{2k}$, the branch
points therefore extend across infinity as a configuration of distinct
points. Over a large circle in the $x$-plane, we can contract this
configuration loop in the $w'$-coordinate to its value at infinity.
Returning to $w=x^{2k}w'$ adds $2k$ full twists. Thus the product of the
local braids is $\tau_{\partial D}^{\,2k}$.
\qed
\section{Boundary framings and braid monodromy}\label{sec:framings}

\subsection{Compatibility and framings}
To read the page monodromy relative to its boundary, we must say how
to follow points on the boundary circles as $x$ moves. In the affine
trivialization of $\mathcal O_{\PP^1}(2k)$, the target coordinate is
$w=s_0/s_\infty$. Choose $R$ so that
$D_x=\{|w|\le R\}$ contains the branch points along the path. The page is
$P_x=\Phi^{-1}(D_x)\cap\pi^{-1}(x)$, with boundary circles covering
$\partial D_x$ in degrees $\lambda_c$. Parametrize the target boundary
circle by $w=Re^{i\theta}$. We obtain the \emph{covering framing} by
lifting this fixed parametrization. Hironaka's framing instead fixes
the angle in the original $y$-plane, $\arg(y-L_i(x))$ at an inner
circle and $\arg y$ at the outer circle \cite[\S2.2]{Hironaka}.
For the corresponding difference between Hironaka's framing and the
framing from a blown-up disk arrangement, see
\cite[Remark~3.6]{PlamenevskayaStarkston}.

By Lemma~\ref{lem:linsys}, $s_0$ vanishes nowhere on the pole
sections $S_c$. Near $S_c$, choose a fiber coordinate $t_c$ vanishing
on $S_c$.  Then $w=c_c\,t_c^{-\lambda_c}(1+O(t_c))$, with a nowhere-zero
leading coefficient $c_c$.

\begin{lemma}\label{lem:framings}
In the covering framing, the page monodromy of $\pi$ along Hironaka's loop
around $x(p)$ is the twist $\tau_{a_p}$, and the framing of the outer circle
$\partial_0$ agrees with Hironaka's.
\end{lemma}

\begin{proof}
Near $\tL_i$ Hironaka uses the fiber coordinate $y-L_i(x)$.  Over $\C$, the
function
\[
 t_i=\frac{y-L_i(x)}{\prod_{p\in\cI\cap L_i}(x-x(p))}
\]
is a fiber coordinate near $\tL_i$ in $\hX$, vanishing on that section: near a point
$p\in\cI\cap L_i$, with $x-x(p)=u$ and $y-L_i(x)=u\,t$, it is a unit
multiple of the blow-up coordinate $t$. Substituting $t_i$ in the expansion
of $w$ gives the leading coefficient
in Hironaka's coordinate
\[
 c_i^{\mathrm H}(x)=c_i(x)\prod_{p\in\cI\cap L_i}(x-x(p))^{\lambda_i}.
\]
Since $c_i$ is nowhere zero on $\C$, its winding number along each loop is
zero. The winding number of $c_i^{\mathrm H}$ is therefore $\lambda_i$
times the sum of the winding numbers about $x(p)$ for $p\in\cI\cap L_i$.
At fixed $\arg w$, the leading term gives a change in $\arg(y-L_i(x))$
equal to the change in $\arg c_i^{\mathrm H}$ divided by $\lambda_i$.
These winding numbers also hold for the boundary circles at sufficiently
large $R$, by homotopy to the leading term. Thus the lifts fix the boundary
pointwise. The loop about $x(p)$ contributes one relative turn on each
$\partial_i$ with $i\in S_p$ and zero on the other inner circles. The
two local monodromies differ by $\prod_{i\in S_p}\tau_{\partial_i}^{\epsilon}$,
with a sign $\epsilon=\pm1$ fixed by the orientation conventions and
independent of $p$.

Near $E_O$ there are no blow-ups at points of $\cI$, and
$w=c_0(x)\,y^{2k}(1+\cdots)$, where $c_0$ is nowhere zero on $\C$.
The same computation shows that the covering framing of
$\partial_0$ agrees with Hironaka's along every loop in $\C$.

By \cite[Lemmas~2.1 and~2.2]{Hironaka}, Hironaka's local monodromy along the
loop around $x(p)$ is $(\prod_{i\in S_p}\tau_{\partial_i})^{-1}\tau_{a_p}$,
and the product along all loops is $\tau_{\partial_0}\prod_i
\tau_{\partial_i}^{-1}$.  In the covering framing the product is the lift of
$\tau_{\partial D}^{2k}$, which is $\tau_{\partial_0}\prod_i
\tau_{\partial_i}^{\mu_i-1}$ by Lemma~\ref{lem:fulltwist}.  The boundary twists
of $\Sigma_{0,s+1}$ generate a free abelian group, and comparing exponents of
$\tau_{\partial_i}$ gives $\epsilon\mu_i-1=\mu_i-1$, so $\epsilon=1$ and the
local monodromy in the covering framing is $\tau_{a_p}$.
\end{proof}

\subsection{Collapsing the inner boundary circles}\label{rem:holes}
For a planar page $P=\Sigma_{0,r}$ with a fixed outer circle $\partial_0$,
let $q$ be the quotient map that collapses the inner boundary circles to
marked points. It induces a homomorphism
\[
 q_*\colon\Mod(P\ \mathrm{rel}\ \partial P)\longrightarrow PB_{r-1}
\]
that kills the inner boundary twists. The resulting braid monodromy is
independent of the inner boundary framings.
For the covering in Theorem~\ref{thm:arrangement}, collapsing the inner
boundary circles of the page over $x$ gives a disk that can be identified
with $\{|y|\le R'\}$ in the vertical line, with the $s$ points $L_i(x)$ marked.

\begin{lemma}\label{lem:collapsed-monodromy}
The image under $q_*$ of the lifted factorization of $w$ in
Theorem~\ref{thm:arrangement} is a braid monodromy factorization of $\cL$
with respect to $\rho$.
\end{lemma}

\begin{proof}
The factors of exponent two and three have identity lift.  Omitting their
loops from the geometric basis defining $w$ leaves a geometric basis of
$\C\setminus x(\cI)$.  By Lemma~\ref{lem:framings}, the outer framing agrees
with that of the coordinate $y$.  After collapsing the inner boundary circles, 
the monodromy along each remaining loop is therefore the braid traced by the points $L_i(x)$ in
$\{|y|\le R'\}$, which is the corresponding braid monodromy factor of $\cL$.
\end{proof}

\section{A Zariski pair of liftable factorizations}\label{sec:zariski}

\subsection{Braid monodromy and the statement}
For the projection $\rho(x,y)=x$ of Section~\ref{sec:arrangement}, the braid
monodromy factorization of $\cL$ records the motion of the points $L_i(x)$
along a geometric basis of loops around $x(\cI)$. We compare these
factorizations up to Hurwitz moves and simultaneous conjugation in $B_s$.

\begin{theorem}\label{thm:zariskiplays}
Let $\cL^+$ and $\cL^-$ be the eleven-line arrangements of
\cite[Theorem~5.3]{ACCM}, denoted $\mathcal H^\pm$ there, in real coordinates
satisfying the hypotheses of Section~\ref{sec:arrangement}, and let $k=30$.  Then
Theorem~\ref{thm:arrangement} gives liftable factorizations $w^+$ and $w^-$ of
$\tau_{\partial D}^{60}$ with respect to one simple branched covering of
degree $224$ with $234$ branch points.  Under any identification of their
pages induced by the covering, the lifted factorizations of $w^+$ and
$w^-$ are not Hurwitz equivalent, even up to simultaneous conjugation by
homeomorphisms of the page.  In particular $w^+$ and $w^-$ are not
equivalent, even up to simultaneous conjugation by braids that carry the
covering to an isomorphic one.
\end{theorem}

For both arrangements, the equations in \cite{ACCM} give
\[
 (\mu_i)=(5,4,5,4,5,6,5,5,6,5,7).
\]
The least common multiple of the numbers $\mu_i-1$ is $60$, so the least
compatible value is $k=30$. For this value, $d=60+\sum_i\lambda_i=224$
and $n=d+s-1=234$. The two arrangements have the same $\mu_i$, so their
coverings for this $k$ are isomorphic.

\subsection{Proof of Theorem~\ref{thm:zariskiplays}}
Suppose that, under an identification of the pages induced by the covering,
the lifted factorizations of $w^+$ and $w^-$ are Hurwitz equivalent up to
simultaneous conjugation by a homeomorphism $h$ of the page.  This
identification matches the boundary circles over the unique cycles of
length $\lambda_0=60$.  Conjugation by $h$ carries
the product of one lifted factorization to that of the other, and both
products equal
$\tau_{\partial_0}\prod_i\tau_{\partial_i}^{\mu_i-1}$.  The twists about the
boundary circles generate a free abelian group.  If $h$ reversed orientation,
conjugation by $h$ would carry this product of positive boundary twists to a
product of negative ones, so $h$ preserves orientation and permutes the
boundary twists.  Moreover,
$\partial_0$ is the only boundary circle with exponent one, since $\mu_i-1\ge3$ for
$\cL^\pm$. Hence $h$ preserves $\partial_0$ and descends to a homeomorphism
of the marked disk obtained by collapsing the inner boundary circles.
The two images under $q_*$ are therefore Hurwitz equivalent up to
simultaneous conjugation in $B_{11}$. By Lemma~\ref{lem:collapsed-monodromy},
the braid monodromies of $\cL^+$ and $\cL^-$ are equivalent.
By Carmona's theorem \cite[Theorem~2]{Carmona}, also stated in
\cite[Theorem~3.8]{Cogolludo}, this would give a homeomorphism
$(\PP^2,\cL^+)\to(\PP^2,\cL^-)$, contrary to \cite[Theorem~5.3]{ACCM}.
An equivalence between $w^+$ and $w^-$ would lift to an equivalence of
their lifted factorizations, since an isomorphism of coverings induces a
homeomorphism of pages.\qed

\begin{corollary}\label{cor:trivialmonodromy}
There is a simple branched covering of the disk of degree $224$ with $446$
branch points and trivial boundary monodromy, whose page is a sphere with
$224$ boundary components. For this covering, the $60$th power of the twist
about the boundary of a proper subdisk admits two liftable factorizations
whose lifted
factorizations are not Hurwitz equivalent, even up to simultaneous
conjugation by homeomorphisms of the page.
\end{corollary}

\begin{proof}
Let $p\colon P\to D$ be the covering of Theorem~\ref{thm:zariskiplays}, with
boundary monodromy $\sigma$ of cycle lengths $\lambda_c$. Embed $D$ in a
larger disk $D'$. In the intervening annulus we extend the cover separately
on the sheets of each cycle. The cycle corresponding to $\partial_c$ is canceled by
$\lambda_c-1$ new branch points, with transpositions whose product is
its inverse. The piece attached to $\partial_c$ is connected, with one
circle covering the inner boundary of the annulus with degree $\lambda_c$
and $\lambda_c$ circles each mapping homeomorphically onto the outer
boundary of the annulus. It has at least eleven boundary circles, since
$\lambda_c\ge10$, and Euler characteristic $1-\lambda_c$ by
Riemann--Hurwitz, so it is planar. Gluing these pieces to $P$ gives a
planar page $P'$ with $224$ boundary circles. The resulting cover
$p'\colon P'\to D'$ has $212$ new branch points, for a total of $446$,
and trivial boundary monodromy. Viewed in the braid group of $D'$, the words $w^\pm$
are liftable factorizations of the same braid with
the same exponents, and their lifted factorizations are the images of those
for $p$ under the homomorphism
$\eta\colon\Mod(P\ \mathrm{rel}\ \partial P)\to
\Mod(P'\ \mathrm{rel}\ \partial P')$ extending by the identity.

Let $c_0,\ldots,c_{11}\subset P'$ be the curves $\partial_0,\ldots,
\partial_{11}$. The $c_c$ are essential and pairwise nonisotopic, and
$\eta$ is injective. Cutting along the $c_c$ separates $P$ from the
attached pieces. The original page $P$ is the only piece with no boundary circle in
$\partial P'$.
Suppose the lifted factorizations of $w^\pm$ for
$p'$ are Hurwitz equivalent up to simultaneous conjugation by a homeomorphism
$h$ of $P'$.  Then $h$ conjugates the product
$\tau_{c_0}\prod_i\tau_{c_i}^{\mu_i-1}$ to itself.  As in the proof of
Theorem~\ref{thm:zariskiplays}, $h$ preserves orientation, and by the
uniqueness of multitwists it permutes the curves $c_c$.
It follows that $h$ preserves $P$ up to isotopy.  The Hurwitz moves act on Dehn twists
about curves in $P$, so the equivalence takes place in the image of $\eta$,
and by injectivity the lifted factorizations of $w^\pm$ for $p$ are Hurwitz
equivalent up to conjugation by $h|_P$, contradicting
Theorem~\ref{thm:zariskiplays}.
\end{proof}

\subsection{Milnor fibers}

Let $W^\pm$ be the total spaces of the planar Lefschetz fibrations
defined by the lifted factorizations in Theorem~\ref{thm:zariskiplays}.
The vanishing cycle $a_p$ represents the nonzero class
$\sum_{i\in S_p}[\partial_i]$ in $H_1(P;\Z)$, so the fibrations
are allowable and give Stein structures on $W^\pm$. Both induce the
open book $(P,\tau_{\partial_0}\prod_i\tau_{\partial_i}^{\mu_i-1})$.

\begin{proposition}\label{prop:milnor}
With these Stein structures, $W^\pm$ are realized by Milnor fibers of
smoothings of possibly analytically distinct rational surface
singularities with reduced fundamental cycle and the same contact link.
\end{proposition}

\begin{proof}
For either arrangement write $L_i=\{y=a_ix+b_i\}$. To build a picture
deformation, give the $i$-th branch of the pencil $C=\bigcup_i\{y=a_ix\}$
the weight $\mu_i$, which counts the marked points with multiplicity
on that branch. The associated rational singularity has reduced
fundamental cycle. Its resolution graph has central weight $-12$ and
eleven arms of $(-2)$-spheres, with lengths determined by the weights
\cite[\S7.1]{PlamenevskayaStarkston}.

Move the lines by setting $C_s=\bigcup_i\{y=a_ix+sb_i\}$. The
parametrization by $x$ of each branch gives a simultaneous normalization.
For $s\ne0$, this is the original arrangement scaled by $s$. Each pair
of branches meets once and transversely, including at $s=0$, so the sum
of the local $\delta$-invariants is constant. Mark each intersection
on every branch through it. On the $i$-th branch the marks are given by
\[
\prod_{p\in\cI\cap L_i}(x-s\,x(p))=0.
\]
This is finite and flat over the parameter disk because the polynomial
is monic in $x$. Its special fiber has length $\mu_i$ at the origin.
For $s\ne0$ it is reduced and marks each preimage of an intersection
once. These divisors contain the relative total multiplicity scheme.
This is a picture deformation in the sense of
\cite[Definition~4.6]{deJongVanStraten} and
\cite[Definition~2.7]{PlamenevskayaStarkston}.

To obtain the Milnor fiber, blow up the bidisk at the marked points and
remove tubular neighborhoods of the proper transforms of the line
disks \cite[Proposition~5.1]{deJongVanStraten}. Projection to the
$x$-disk gives a Lefschetz fibration compatible with the Stein structure
of the Milnor fiber \cite[Lemma~3.2]{PlamenevskayaStarkston}. Under dilation, its vanishing
cycles are the transported circles $a_p$ of the original arrangement.
The boundary framing difference acts by boundary twists, which fix
the isotopy classes of these cycles. This gives the fibration defining
$W^\pm$.
\end{proof}

\begin{lemma}[A test for simple connectivity]\label{lem:meridian-elimination}
Let $\mathcal A$ be a real affine line arrangement. Its complex lines cut
out disks in a large bidisk in $\C^2$. Blow up their intersection points
to separate the disks, and let $W(\mathcal A)$ be the complement of tubular
neighborhoods of the resulting disjoint disks.

In the real plane, start with all the lines of $\mathcal A$ and mark their
intersection points. Apply the following rules, keeping the marked points
as lines are deleted.
\begin{enumerate}
\item If just one remaining line passes through a marked point, delete it.
\item If a bounded complementary region (a chamber) of the remaining
arrangement has an even number of sides, with three remaining lines
through one vertex and two through each of the others, delete the third
line at that vertex, the one that is neither of the two sides.
\end{enumerate}
The set of lines left when neither rule applies is independent of the
order of deletion. If no lines remain, $W(\mathcal A)$ is simply connected.
\end{lemma}

\begin{proof}
In the four-manifold, deleting a line from the real diagram means
gluing back the neighborhood of its disk. The meridian of the disk then
bounds a normal disk in the restored neighborhood, so by van Kampen the
fundamental group is unchanged if the meridian is already trivial.
At a marked point, the
blow-up replaces the crossing by an exceptional sphere. Each disk still
missing from the complement leaves a hole in that sphere. If only one
disk is missing, the sphere with one hole is a disk that bounds its
meridian. This proves the first rule.

For the second rule, a marked point in the interior of a side lets us
delete that side by the first rule. The first rule then propagates the
deletion around the double vertices to the other sides and finally to
the third line at the triple vertex. Otherwise, choose a base point
inside the region and take paths within it to the sides. Detours around
interior marks do not change the based meridians, since a small blown-up ball
about each such mark is simply connected and misses the remaining disks.
At a corner, the directions into the region form
an interval on the exceptional sphere between the holes for the two sides.
The paths in the corner can be slid onto this interval, so the side
meridians agree with the corresponding based meridians on the sphere.
At a double vertex the sphere with two holes is an annulus, and the two
meridians are inverse. Go around the boundary from one side at the triple
vertex to the other, passing through an odd number of double vertices.
The meridian is inverted at each step, so the two side meridians at the
triple vertex are inverse. In the relation on the sphere with three
holes, they cancel in either cyclic order, leaving the third meridian
trivial. Thus every deletion preserves the group. If all disks are put
back, only the simply connected blow-up of the bidisk remains.

To compare two orders, run the first sequence until neither rule applies.
Consider the first line in the second sequence that survived the first.
All earlier deletions in the second sequence have already been made in
the first. If a marked point allowed this line to be deleted in the
second sequence, the same point allows it at the end of the first. If a
region allowed the deletion and all its sides survive, the same region
still allows it. If a side is gone, the chain of double vertices forces
all the sides to be gone by the first rule. The marked point at the triple
vertex then allows deletion of the third line. In every case a rule
still applies, contradicting termination. Interchanging the sequences
gives the same set of surviving lines.
\end{proof}

\emph{Application to $W^\pm$.}
Use the blown-up disk model of Proposition~\ref{prop:milnor} in the affine
chart $3x+y+z=1$, with the line equations of
\cite[Proposition~3.7 and \S5]{ACCM}. The following quadrilaterals
have the indicated sides in cyclic order. Their vertices are the
successive intersections of the sides. Substituting the vertices in
the line equations shows that each side bounds the quadrilateral and
that no other line crosses its interior.
\[
\begin{array}{c|c|c}
\text{Arrangement} & \text{Sides in cyclic order} &
  \text{Third line at the triple vertex}\\ \hline
\cL^+ & M_1,\ N,\ M_5,\ L_5 & L_4\\
\cL^- & M_1,\ N,\ M_3,\ L_3 & L_2\\
\cL^-\setminus\{L_2\} & M_1,\ L_3,\ M_3,\ L_5 & L_4
\end{array}
\]
The line at infinity of this chart avoids all intersection points.
The disk-arrangement model is unchanged up to isotopy along a path of
generic complex affine charts and projections. The first row deletes
$L_4$ in the plus arrangement. In the minus arrangement the second
row deletes $L_2$, and the third then deletes $L_4$.

Once $L_4$ is gone, the double point $L_4\cap N$ lets us delete $N$.
The double points of $N$ with $M_1,M_3,M_5,L_1$ then delete those
four lines. The double points $M_4\cap L_4$,
$M_3\cap L_3$, and $M_5\cap L_5$ give three more deletions. The triple
points on $M_2,L_1,L_5$ and on $L_2,M_3,L_5$ account for $M_2$ and
$L_2$, the latter already deleted in the minus arrangement. Thus both
$W^+$ and $W^-$ are simply connected.

\begin{corollary}\label{cor:fillings}
The Milnor fibers $W^\pm$ arise from picture deformations with the same
incidence matrix. The boundary identification given by the common cap in
the proof extends to an orientation-preserving homeomorphism $W^+\to W^-$.
The fillings belong to distinct strong symplectic deformation classes.
\end{corollary}

\begin{proof}
Let $Y$ be the common boundary with its planar open book. The Stein
forms on $W^\pm$ are exact, so both fillings are symplectically minimal.
For an ordinary planar open book, the classification theorem of Lisi,
Van Horn-Morris, and Wendl identifies the strong symplectic deformation
classes of minimal fillings with the isomorphism classes of allowable
Lefschetz fibrations inducing that open book \cite[Theorem~1.5]{LVW}.
See also \cite[Theorems~1 and~3]{Wendl} for existence and deformation
stability in the planar case.
An isomorphism of the two fibrations would give Hurwitz equivalence up
to simultaneous conjugation by a page diffeomorphism, contrary to
Theorem~\ref{thm:zariskiplays}.

The rows of the incidence matrix correspond to the disks and the
columns to the marked points. An entry is $1$ if the disk passes through
the point and $0$ otherwise. The two arrangements have the same labeled
incidences, so the matrices agree. Both fillings are simply connected by
the application of Lemma~\ref{lem:meridian-elimination}. To apply Boyer's
theorem, we compare their intersection forms together with the boundary
identification. For this we use the cap $U$ of
\cite[\S6.3]{PlamenevskayaStarkston}, obtained
by attaching one $2$-handle per line to a $4$-ball. Its attaching link
is the boundary link of the line disks. For either arrangement this
link is isotopic to a union of labeled fibers of the Hopf fibration.
The matching framings $-\mu_i$ therefore give the same cap for both
fillings. Gluing $U$ to $W^\pm$ gives
$X^\pm=\#_{p\in\cI}\overline{\mathbb{CP}^{\,2}}$.
The inclusion of $H_2(U)$ in $H_2(X^\pm)$ is given, in the exceptional
bases, by the same incidence data. Matching the exceptional classes
therefore gives a lattice isometry that fixes the image of $H_2(U)$
and restricts to an isometry of the orthogonal complements $H_2(W^\pm)$.
The form on $H_2(U)$ has diagonal entries $-\mu_i$ and off-diagonal
entries $-1$, since two lines meet once. Its negative is the sum of
$\operatorname{diag}(\mu_i-1)$ and the all-ones matrix. This is positive
definite with determinant greater than one, so $Y$ is a rational
homology sphere with nontrivial integral first homology.

A $2$-cycle in $X^\pm$ transverse to $Y$ splits into relative cycles on
the two sides, with opposite boundaries in $Y$. By duality,
the two relative classes are computed by pairing the original class
with $H_2(W^\pm)$ and $H_2(U)$, respectively. Every relative class in
$H_2(W^\pm,Y)$ comes from $H_2(X^\pm)$ since $H_1(U)=0$. The ambient
isometry preserves the cap-side class in the common group $H_2(U,Y)$,
so it also preserves the boundary class on the filling side. Thus the
isometry of the filling forms induces the prescribed boundary map on
$H_2(W^\pm,Y)/H_2(W^\pm)\cong H_1(Y)$. Both
Kirby--Siebenmann invariants vanish because the fillings are smooth.
Boyer's theorem now extends the boundary identification to a
homeomorphism of the fillings
\cite[Theorem~0.7 and Proposition~0.8(i)]{Boyer}.
\end{proof}

\begin{samepage}
\begin{remark}
After taking the same finite number of interior connected sums with
$S^2\times S^2$ on each side, the boundary identification extends to a
diffeomorphism. Indeed, the homeomorphism of Corollary~\ref{cor:fillings}
has vanishing Casson--Sullivan invariant, since
$H^3(W^+,Y;\Z/2)\cong H_1(W^+;\Z/2)=0$. The relative stable smoothing
theorem then applies \cite[Proposition~2.27]{Galvin}.
\end{remark}
\end{samepage}

\begin{remark}
For the singularity link of Proposition~\ref{prop:milnor}, the absolute
central weight is one more than the number of arms. In the compactification of
\cite[\S2.2]{ChoiPark}, blowing down the exceptional spheres gives the
original line arrangement together with the line at infinity. The added
line meets the eleven lines at distinct points. The original
arrangement already has more than one point of multiplicity at least three.
Thus \cite[Theorem~1.1 and Remark~1.2]{ChoiPark} rules out obtaining either
filling from the minimal resolution by a sequence of symplectic rational
blowdowns.
\end{remark}

\section*{Declaration on the use of AI tools}
The author used Claude (Anthropic) and Codex (OpenAI) to assist with
literature searches, checks of computations and arguments, and editing.
The author takes full responsibility for the mathematical claims and the
final text.

\end{document}